\documentclass[11pt, a4paper]{article}

\usepackage{amsmath,amssymb, amsthm}
\usepackage{inputenc}
\usepackage{mathrsfs}
\usepackage{fullpage}
\usepackage[affil-it]{authblk}

\usepackage{graphicx}
\graphicspath{ {images/} }

\usepackage{wrapfig}

\usepackage{dsfont}
\usepackage{tikz}

\usetikzlibrary{arrows,calc}

\usepackage{verbatim}

\usepackage{hyperref}

\usepackage{xcolor}
\hypersetup{
  colorlinks   = true, 
  urlcolor     = {blue!90!black}, 
  linkcolor    = {blue!90!black}, 
  citecolor   = {red!90!black} 
}

\usepackage{setspace}
\newcommand{\R}{\mathbb{R}}
\newcommand{\N}{\mathbb{N}}
\newcommand{\Z}{\mathbb{Z}}
\newcommand{\Go}{\Gamma _0}

\renewcommand{\d}{\mathrm{d}} 
\newcommand{\B}{\mathbf{B}}

\theoremstyle{definition}
\newtheorem{defi}{Definition}[section]

\theoremstyle{plain}
\newtheorem{thm}[defi]{Theorem}
\newtheorem{prop}[defi]{Proposition}
\newtheorem{lem}[defi]{Lemma}

\theoremstyle{definition}

\newtheorem{con}[defi]{Condition}
\newtheorem{rmk}[defi]{Remark}

\author{Viktor Bezborodov \thanks{Email: \texttt{viktor.bezborodov@univr.it}}}

\author{
Luca Di Persio \thanks{Email: \texttt{luca.dipersio@univr.it}}}

\affil{
\emph{University of Verona - Department of Computer Science}}

\title{Maximal irreducibility measure
\\ for spatial birth-and-death processes}

\date{}
\begin{document}

\maketitle
\noindent
{\small
}

\begin{abstract}

We prove that
a spatial birth-and-death process is both $\phi$-irreducible and $\psi$-irreducible
under rather general conditions on the birth and death rates.
It is also shown that
every maximal irreducibility measure is equivalent to 
the Lebesgue-Poisson measure on the space of finite configuration.

\end{abstract}

\section{Introduction}

The basic question of stochastic stability analysis
for a Markov process is whether the chain is irreducible.
The notion of irreducibility for countable state space
Markov processes is not directly transferable 
to Markov processes with continuous state spaces. 
The most widely used generalization 
is the so called $\phi$-irreducibility,
see, e.g., Myen and Tweedie \cite{MeynTweedieBook}.
The aim of this paper is to prove that under certain general conditions
the Lebesgue-Poisson measure
is a maximal irreducible measure for
continuous-space
birth-and-death processes.
Roughly speaking it means that, whatever the initial condition is,
a set will be hit by the process with positive probability 
if and only if it is of positive Lebesgue-Poison measure.

We describe and define spatial birth-and-death processes 
in Section \ref{laxative}.
The pioneering works on spatial birth-and-death processes are 
Preston \cite{Preston} and 
Holley and Stroock \cite{HolleyStroock}.
More recent studies of various related aspects include for example
\cite{FournierMeleard, GarciaKurtz, BaDdynamics}, while we refer the interested reader to \cite{DiPersioRW}, for
questions related to the connections with the theory of {\it random walks in randome media}, and to \cite{Barbu, Benazzoli}, and references therein, 
to what concerns the links with more applicative problems as those arising in financial and neurobiological settings.

The paper is organized as follows:
in Section \ref{splatter}
we recall the notions of 
$\phi$-irreducibility and maximal irreducibility for
measures; in Section \ref{dibber}
we recall the definition of the Lebesgue-Poisson measure; in Section \ref{laxative} we describe 
the birth-and-death processes we consider 
and give our main result, Theorem \ref{irreducibility thm}; the proofs
are collected  in  Section \ref{proofs}.

\section{Irreducible and maximal irreducible measures}\label{splatter}
In what follows we shall adopt the notation used in 
\cite{MeynTweedieBook}.
Let $X$ be a Polish space and $\mathscr{B} (X)$ 
be its Borel $\sigma$-algebra. 
 We will consider a
 Markov chain with transition probability kernel $P$
 and initial distribution $\mu$  defined on the 
 canonical space $\Omega = \prod _{i=0} ^\infty X $,
 with $\Phi _n$ being the coordinate mappings,
 \[
 \Phi _n ((x_0,x_1,...)) = x_n.
 \]
The corresponding measure will be denoted by $P _{\mu}$,
so that
 for any Borel sets $A_0,...,A_n \in \mathscr{B} (\Omega)$,
\begin{align}
 \begin{aligned}
 P_\mu (\Phi  _0 \in & A _0 , \Phi  _1 \in A _1,...,\Phi  _n \in A _n)= \\
  & \int _{y_0 \in A_0} ... \int _{y_{n} \in A_{n}} \mu(dy_0)P(y_0, dy_1)...P(y_{n-1}, dy_n).
 \end{aligned}
 \end{align}
 Let $P_x$ denote the distribution of $\Phi $ in $\Omega$ 
 when the initial distribution is the Dirac measure at $x$,
 $P_x \{ \Phi  _0 = x  \} =1$.  For any set $A \in \mathscr{B} (X)$,  
$
 \tau _A = \min \{ n \geq 1: \Phi  _n \in A \}  
 $
 is
 called the \textit{first return time}.
 Define also the return probabilities

 \begin{align}
  \begin{aligned} L(x,A) :&=   P_x \{ \tau _A <\infty \}   
                          = P_x \{ \Phi  \text{ ever enters } A  \}.
                                         \end{aligned}
 \end{align}

 \begin{defi}\label{MeasuresDefinition}
 A finite non-trivial measure $\phi$ is
 called 
 $\phi$-\textit{irreducible} for the chain $\Phi$ 
 if $\phi (A) >0$ implies that  
 $$
 L(x, A)>0,  \ \ \ x \in X.
 $$ 
  A finite non-trivial measure $\psi$ is
 called  
 $\psi$-\textit{maximal irreducible} for the chain $\Phi$
 if
 $$ 
 (\forall x \in X : \ L(x, A)>0 ) \Leftrightarrow \psi (A) > 0.
 $$
 
 \end{defi}
 
 The measures $\phi$ and $\psi$ from the above definition
 are called an \textit{irreducibility measure} and
 a \textit{maximal irreducibility measure} for $\Phi $, respectively.
The next proposition provides a sufficient condition for an irreducibility
measure to be a maximal irreducibility measure.

  \begin{prop} \label{lush}
   If $\Phi $ is $\phi$-irreducible and the measure $\phi$ is such that
   $\phi \{ y: P(y,A) >0 \} = 0$
   whenever $\phi (A) = 0$,
   then $\Phi $ is $\psi$-irreducible with $\psi = \phi $. 
   
  \end{prop}

\section{Lebesgue-Poisson measure}
\label{dibber}

The state space of a continuous-time, continuous-space
birth and death process is
 \[
 \Gamma _0(\R ^\d)=\{ \eta \subset \R ^\d : |\eta| < \infty \},
\]
where $|\eta|$ is the number of points of $\eta$.
$\Gamma _0(\R ^\d)$ is often called the 
\emph{space of finite configurations}.
The space of $n$-point configuration is
$\Gamma _0 ^{(n)} (\R ^\d) :=
\{ \eta \subset \R ^\d : |\eta| =n \}
\subset \Gamma _0(\R ^\d)$.
We will use $\Gamma _0$ and $\Gamma _0 ^{(n)}$ as  shorthands
for $\Gamma _0(\R ^\d)$ and $\Gamma _0 ^{(n)}(\R ^\d)$,
respectively.
  For  $\eta, \zeta \in \Go$, $|\eta|=|\zeta|>0$, 
we define
\begin{equation}\label{chirp}
  \rho(\eta,\zeta) :=   
  \min\limits _{\varsigma} \max\limits _{x\in \eta}
  \{|\varsigma(x) - x| \},
\end{equation}
where 
minimum is taken over the set of all bijections 
$\varsigma : \eta \to \zeta$. Note than 
in \eqref{chirp} the notation $|\cdot|$
is used for the Euclidean distance in $\R ^\d$
(as opposed to the number of points as in $|\eta|$),
which hopefully should not lead to ambiguity.
For  $\eta \in \Go$ and $a>0$, let 
\[
 \mathbf{B}_{\rho} (\eta , a): = \{ \zeta \in \Gamma ^{(|\eta|)}_0 \mid \rho(\eta,\zeta) \leq a \}.
\]

The $\sigma$-algebra 
can be defined as
\[
 \mathscr{B}(\Gamma _0 ) = 
 \sigma \left( \{\varnothing \}, 
 \mathbf{B}_{\rho} (\eta , a), \eta \in \Go, a>0
  \right).
\]

Let
\begin{equation}
\widetilde {(\R^\d)^n} := \{ (x_1 ,..., x_n) \in (\R^\d)^n \mid x_j \in \R ^\d, j=1,...,n, 
x_i \ne x_j , i \ne j \},
\end{equation}
and let $sym$ be the mapping
\begin{equation*}
  \begin{split}
   \bigsqcup\limits _{n=0}^\infty \widetilde {(\R^\d)^n} 
   \ni (x_1,...,x_n) \mapsto \{x_1 ,..., x_n  \} \in \Go. \\
  \end{split}
 \end{equation*}
 We are now going to  define the Lebesgue-Poisson measure on $(\Go, \mathscr{B}(\Go))$. 
 For any $n \in \N$, let
 $l ^{\otimes n} _\d $ be the restriction 
 of the Lebesgue measure to $\widetilde {(\R^\d)^n}$.
 We denote
 by 
 $\lambda ^{ (n)} $
 the projection of this measure on 
 $\Gamma _0 ^{(n)}$,
 
 \[
  \lambda ^{ (n)} (A) = l _\d ^{\otimes n} (sym^{-1} A) \;,\;
 A \in \mathscr{B}(\Gamma _0 ^{(n)}).
 \]

 On $\Gamma _0 ^{(0)}$
 the measure $\lambda ^{ (0)} $ is given by 
 $\lambda ^{ (0)}(\{ \varnothing \}) =1$.
 The \textit{Lebesgue-Poisson measure} on $(\Go, \mathscr{B}(\Go))$
 is defined as 
 
 \begin{equation}
  \lambda  := \sum\limits _{n=0} ^\infty \frac{1}{n!} \lambda ^{ (n)}.
  \end{equation}
 Let us note that the
measure $\lambda$ is infinite.

\section{Birth-and-death processes and the main result}
\label{laxative}

Denote by $\mathscr{B}(\R ^\d)$ the Borel $\sigma$-algebra on $\R ^\d$.
 The evolution of a spatial birth-and-death process admits the following
description.
Two functions characterize the 
 development in time,
 the birth rate 
 $b: \R ^\d  \times \Gamma _0 (\R^\d) \rightarrow [0,\infty)$
 and the death rate 
 $d: \R ^\d \times \Gamma _0 (\R^\d) \rightarrow [0,\infty)$. 
 If the system is in state $\eta \in \Go$ 
 at time $t$, then the probability 
 that a new particle is added (``birth'' event)
 in a bounded set $B\in \mathscr{B}(\R ^\d)$
 over the time interval $[t;t+ \Delta t]$ is 
 \[
   \Delta t \int\limits _{ B}b(x, \eta)dx + o(\Delta t),
 \]
 while the probability that
a particle $x \in \eta$ is removed from the configuration ( ``death'' event), over 
 time interval $[t;t+ \Delta t]$ is
 \[
  d(x, \eta) \Delta t + o(\Delta t),
 \]
and simultaneous events cannot occur.
In other words,
the rate at which a birth 
 occurs in $B$ is $\int_B b(x, \eta)dx$,
 and the rate at which a particle $x\in \eta$
dies is $d(x,\eta)$, 
and no two events happen at the same time.
Various aspects of birth-and-death processes 
are considered in, e.g., \cite{FournierMeleard, GarciaKurtz, KondSkor}.
Here we focus our attention on the  embedded Markov chain
of the birth-and-death process, namely the Markov chain on $\Go$
with transition probabilities
 \begin{equation}\label{transition probabilities}
 \begin{gathered}
    Q(\eta, \{ \eta \setminus \{ x \} \}) = 
  \frac{d(x,\eta)}{(B+D)(\eta)} , \ \ \ \ \ \ \ x \in \eta, \ \ \eta \in \Go,  \\
  Q(\eta, \{ \eta \cup \{ x \}, x \in U\} )
  = \frac{\int _{x\in U} b(x, \eta) dx}{(B+D)(\eta)}, \ \ \
  U \in \mathscr{B}(\R ^\d), \eta \in \Go, 
 \end{gathered}
 \end{equation}
 where $(B+D)(\eta) = 
 \int\limits _{x \in \R^\d} b(x, \eta) dx + \sum\limits _{x \in \eta} d(x, \eta ) $
 is the jump rate at $\eta $.

 Denote by $Q _{\alpha}$
    the distribution 
    of the Markov chain on
    ${((\Go)^\infty, \mathcal{B}((\Go)^\infty) )}$
    with transition probabilities 
    \eqref{transition probabilities}
    and initial value $\alpha \in \Go$.
     Here $\mathcal{B}((\Go)^\infty) )$ is the $\sigma$-algebra
    generated by the coordinate mappings. 
    Let $(\xi _n)_{n\in \Z_+}$
    be the coordinate mappings
    $((\Go)^\infty, \mathcal{B}((\Go)^\infty) )$, that is
    $
     \xi _n (\boldsymbol{\eta}) = \eta _n
    $ for $\boldsymbol{\eta} = (\eta _0, \eta _1, ...) \in (\Go)^\infty$.
     Under $Q _{\alpha}$,
     $( \xi _n  )_{n \in \Z _+}$
     is a Markov chain with transition probabilities
     \eqref{transition probabilities}.

Concerning the functions $b$ and $d$, we assume that they are  continuous functions in both variables, satisfying the following conditions

\begin{con}[Sublinear growth] There 
 exist $c_1, c_2 >0$ such that

\begin{equation} \label{sublinear growth for b} 
\int\limits _{\R ^\d} {b}(x, \eta ) dx \leq c_1|\eta| +c_2.
\end{equation}
\end{con}
\begin{con}
 
We require

\begin{equation} \label{condition on d}
\forall m \in \N : \sup\limits _{x\in \R^\d,|\eta| \leq m} d(x, \eta) < \infty.
\end{equation}
 \end{con}
 \begin{con}[Non-degeneracy of $d$] The infimum
   \begin{equation} \label{inf d > 0}
    \inf\limits _{\eta \in \Gamma _0 (\R^\d), x \in \eta} d(x,\eta) >0,
  \end{equation}
   \end{con}
    \begin{con}[Non-degeneracy of $b$]
  For some constants $r>0$ and $c_3 >0$,
\begin{equation} \label{assump on b}
\begin{split}
   b(x, \eta )   > c_3  , 
   \text{if there exists \ } y \in \eta, |x-y| \leq r,  \\
   \text{ and \ }
   b(x, \varnothing )   > c_3 \text{ for } 
   x \in B_{\varnothing}, B_{\varnothing}
   \text{ is some open ball in } R^\d.
   \end{split}
  \end{equation}
  \end{con}
The following theorem constitutes the main result of the present paper.
  \begin{thm}\label{irreducibility thm}
 
 The Lebesgue-Poisson measure $\lambda$
 is a maximal irreducibility measure for $(\xi _n)_{n\in \N}$.

 \end{thm}

 In other words,
    
    \[
  (\forall \alpha :  Q _{\alpha} \{ (\xi _n)_{n\in \Z_+} \text{ ever enters } A \} >0 ) 
  \Leftrightarrow \lambda (A) >0.
    \]

  \begin{rmk}
  The second part of  \eqref{assump on b}
   means that points may come ``out of nowhere''.
   We need such kind of condition in order 
   for $\varnothing$ not to be an absorbing state of
   the Markov chain $(\xi _n)_{n \in \N}$. 
   Also, each of conditions 
   \eqref{inf d > 0} and
   \eqref{assump on b} implies that
   every state $\eta \in \Go$, $\eta \ne \varnothing$,
   is non-absorbing.

  \end{rmk}

\section{Proofs}
\label{proofs}

 \textbf{Proof of Proposition \ref{lush}}. Let $\phi$ be a measure
  satisfying conditions of the proposition. We first prove that 
  
  \begin{equation} \label{mundate}
  \phi \{ y: L(y,A) >0 \} = 0 \text{ whenever }\phi (A) = 0.
  \end{equation}
  Note that 
  
  \begin{equation} 
  \{ y: L(y,A) >0 \} = \bigcup\limits _{n \in \N } \{ y:  P^n(y,A)  >0 \}.
  \end{equation}
  
  For $A \in \mathscr {B} (X)$ and $k \in \N$, denote $A^{(-k)} := \{x\in X : P^k(x,A) >0 \}$. To
  prove \eqref{mundate}, we will proceed 
  by induction and show that $\phi \{ y: P^n(y,A) >0 \} = 0$ as long as $\phi (A) = 0$,
  for all $n \in \N$.
  Assume that $\phi \{ y: P^m(y,A) >0 \} = 0$ whenever $\phi (A) = 0$. Then, if $\phi (A) = 0$,
  
  $$
  \phi \{ y: P^{m+1}(y,A) >0 \} = \phi \{ y: \int\limits _{x \in X} P(y,dx)P^m(x,A)  >0 \} \leq
  $$ 
  $$
  \phi \{ y: \int\limits _{x \in X} P(y,dx)I_{A^{(-m)}}(x)  >0 \} = 
  \phi \{ y:  P(y,A^{(-m)}) >0 \} =0.
  $$
 
 The base case is given in the condition, therefore \eqref{mundate} holds. 
 
 Assume now that the statement of the proposition does not hold, so that $\phi$ 
 is not a maximal irreducible 
 measure for $\Phi $. Proposition 4.2.2 from
 \cite{MeynTweedieBook} 
 implies the existence of a maximal irreducible measure $\psi'$ for $\Phi $. 
 Then there exists 
 a set $C \in \mathscr{B} (X)$ such that $\phi (C) = 0$ whereas $\psi' (C) >0$.
 By definition of irreducibility,
 $L(x,C)>0$ for all $x \in X$. By \eqref{mundate},
 $\phi \{ y: L(y,C) >0 \} = 0$, hence $\phi (X) = 0$, which contradicts
 to the non-triviality of $\phi$.
 $\qed$

Define a
 \textit{path} of configurations as a finite sequence 
 of configurations $\zeta _0, \zeta _1, ... , \zeta _n$
 such that $|\zeta _k \bigtriangleup \zeta _{k+1}| =1$,
 $k=0,...,n-1$, and 
 if $\zeta _{k+1}= \zeta _k \cup \{z\}$, then 
 $|z - y | \leq \frac r2$ for some $y \in \zeta _k$;
 that is, $\zeta _{k+1}$ is obtained from $\zeta _k$ either by
 adding one point to $\zeta _k$ or by removing one
 point from $\zeta _k$; in the case of the adding, it is required
 that the ``new'' point appears not further than $\frac r2$
 from an ``old'' one. If $\zeta _k =\varnothing$, then we require
 $\zeta _{k+1} = \{ x_{\varnothing} \}$, where $x_{\varnothing}$
 is the center of $ B_{\varnothing}$.
 We say that such a path has length $n$, 
 and we call $\zeta _0 $ and $ \zeta _n$ the starting vertex
 and the final vertex, respectively. Also, we say that 
 $\zeta _0, \zeta _1, ... , \zeta _n$ is a path 
 from $\zeta _0 $ to $ \zeta _n$.

\begin{lem}\label{5lem11}
 For all $\eta \in \Go$ there exists a path
 from $\varnothing$ to $\eta$.
 \end{lem}
 
  \textbf{Proof}. We will show that there exists a path 
  from $\varnothing$ to $\eta$
 of length less than
 \[
  2\bigl(\sum\limits _{x \in \eta } |x - x_{\varnothing}|\frac 4r  +|\eta| \bigr),
 \]
where $x_{\varnothing}$ is the center of $B_{\varnothing}$.

  Starting from $\varnothing$ and only adding points,
  we see that there 
  exists a path of length 
 \[
 \leq \bigl(\sum\limits _{x \in \eta } |x - x_{\varnothing}|\frac 4r 
 +|\eta| \bigr),
 \]
 with the starting vertex $\varnothing$
 and with the final vertex being some
 configuration $\eta ' \supset \eta$. Indeed, for each $x \in \eta$
 there exists a sequence of points $x_{\varnothing} = x_0, x_1,...,x_n =x$
 such that $|x_i - x_{i+1}| \leq \frac r4$ and 
 $n \leq |x - x_{\varnothing}|\frac 4r$.
 Having reached $\eta ' \supset \eta$, we only need 
 to delete some points from $\eta '$. $\qed$

\begin{lem}\label{5lem12}
 Let $\varnothing = \eta _0, \eta _1, ... , \eta _n$ be a path.
 Then for every $a >0$
 \[
  Q^n (\eta _0, \B_{\rho} (\eta_n , a)) > 0.
 \]
\end{lem}
\textbf{Proof}.
Without loss of generality we can assume $a< \frac r4$.
Denote $ A_k = \B_{\rho} (\eta_k , a )$. 
We will first show that

\begin{equation}\label{tick}
  \inf _{\eta \in A_k} Q( \eta , A_{k+1}) \geq \bar c _n
\end{equation}
for some positive constant $\bar c _n$ that 
depends  on $n$ but does not depend on the path
we consider. 

We have either $\eta _{k}\subset \eta _{k+1}$ or $\eta _{k}\supset \eta _{k+1}$.
Consider first the case $\eta _{k}\subset \eta _{k+1}$.
We know that $\eta _{k+1} = \eta _{k} \cup \{z\}$, where 
$|z-y|\leq \frac r2$
for some $y \in \eta _k$.

Take arbitrary $\eta \in A_k$. There exists $y' \in \eta$
such that $|y - y'| \leq a$. For $x\in B_{a}(z)$
we have then $|x - y'| \leq |x - z|+|z - y|+|y - y'|
\leq a + \frac r2 +a <r $. Moreover, if $x\in B_{a}(z)\setminus \eta$,
then $\eta \cup \{ x \} \in A_{k+1}$. 

From \eqref{assump on b} we obtain
\[
 Q( \eta , A_{k+1}) \geq \frac{\int\limits _{x\in B_{a}(z)} b(x, \eta) dx}{(B+D)(\eta)}
 \geq \frac{\int\limits _{x\in B_{a}(z)} c_3 dx}{(B+D)(\eta)} 
\]
\[
 = \frac{c_3 a^\d v_\d}{(B+D)(\eta)} ,
\]
where $v_\d$ is the volume of a unit ball in $\R ^\d$. By
\eqref{condition on d},
the denominator of the last fraction is bounded in $\eta$,
$\eta \in \bigsqcup _{k=0}^{n} \Gamma ^{(k)}_0(\R ^\d)$.
Therefore, \eqref{tick} holds.

Now we turn our attention to the case when
 $\eta _{k} \supset \eta _{k+1}$.
We may write 
$\eta _{k+1} = \eta _{k} \setminus \{ y \}$
for some $y \in \eta _{k}$, and 
\eqref{tick} follows from \eqref{inf d > 0}.

The statement of the lemma follows from \eqref{tick}, 
since

\[
Q ^{n}(\varnothing, \B_{\rho} (\eta_n , a)) =   \int\limits _{\zeta_1, \zeta _2, ... , \zeta _{n}} 
 {Q(\varnothing, d \zeta _1)}   
 {Q(\zeta _1, d \zeta _2)}
 {Q(\zeta _2, d \zeta _3)} \times ...
 \times
{Q(\zeta _{n-1}, d \zeta _{n})} I_{\{ \zeta_{n} \in \B_{\rho} (\eta_n , a) \}} 
\]
\[
 \geq \int\limits _{\zeta_1, \zeta _2, ... , \zeta _{n}} 
 {Q(\varnothing, d \zeta _1)}   
 {Q(\zeta _1, d \zeta _2)}
 {Q(\zeta _2, d \zeta _3)} \times ...
 \times
{Q(\zeta _{n-1}, d \zeta _{n})} I_{\{ \zeta_{k} \in \B_{\rho} (\eta_k , a), k=1,...,n \}} \geq 
(\bar c_n )^n.
\]

\begin{lem}\label{5lem13}
Let $A \in \mathscr{B}(\Go)$, 
  $\beta ' \in \Gamma _0 ^{(n)}$ and 
  $\lambda (A \cap \B _{\rho} (\beta ' , \frac r4)) >0$.
  Then
  \[
   Q^{2n}(\beta, A) >0
  \]
 for any $\beta \in \B _{\rho} (\beta ' , \frac r4) $.
  
\end{lem}

\textit{The idea of the proof}. Let $\beta = \{x_1,...,x_n \}$.
The event $R$ described in the next sentence has
positive probability.
Let 
$\xi _1 = \beta \cup \{y_1\}$
for some $y _1 \in B_{\frac r4}(x_1)$,
$\xi _2 = \xi _1 \setminus \{x_1\}$,
$\xi _3 = \xi_2 \cup \{y_2\}$
for some $y _2 \in B_{\frac r4}(x_2)$,
$\xi _4 = \xi _3 \setminus \{x_2\}$,
and so on, so that
$\xi _{2n} = \xi _{2n-1} \setminus \{x_{n}\} $.
We will see that 
$Q_{\beta}\{\xi _{2n} \in A \mid R \} >0$.

\textbf{Proof}.
Fix $\beta = \{x_1,...,x_n \}$.
Consider a measurable subset $\Xi$ of 
 $(\Gamma _0)^{(2n)} $,
 \begin{align*}
  \Xi = \bigg\{ (\zeta_1,...,\zeta_{2n}) \mid 
  \zeta _{2k-1} = \{y_1,...,y_k,x_k,...,x_n \},
  \zeta _{2k} = \{y_1,...,y_k,x_{k+1},...,x_n \}, k=1,...,n, \\
  \text{ for some distinct } y_1,..., y_n
  \in \R ^\d \text{ satisfying }  |y_k - x_k|\leq \frac r4
  \bigg\}.
 \end{align*}
 Define $R = \{(\xi _1,...,\xi _{2n}) \in \Xi \}$.

 By the Markov property,
 \begin{equation}\label{infidelity}
 \begin{gathered}
 Q ^{2n}(\beta, A) =   \int\limits _{\zeta_1, \zeta _2, ... , \zeta _{2n}} 
 {Q(\beta, d \zeta _1)}   
 {Q(\zeta _1, d \zeta _2)}
 {Q(\zeta _2, d \zeta _3)} \times ...
 \times
{Q(\zeta _{2n-1}, d \zeta _{2n})} I_{\{ \zeta_{2n} \in A \}}
\\
\geq \int\limits _{\zeta_1, \zeta _2, ... , \zeta _{2n}} 
 {Q(\beta, d \zeta _1)}   
 {Q(\zeta _1, d \zeta _2)}
 {Q(\zeta _2, d \zeta _3)} \times ...
 \\
 \times
{Q(\zeta _{2n-1}, d \zeta _{2n})} I_{\{(\zeta_1,...,\zeta _{2n}) \in \Xi  \}} 
I_{\{ 
(\zeta _1 \setminus \beta) \curlyvee (\zeta _3 \setminus \zeta _2) \curlyvee
... \curlyvee (\zeta _{2n-1} \setminus \zeta _{2n-2})
\in  sym ^{-1} A \}}. 
  \end{gathered}
  \end{equation}
  Here for singletons $\mathbb{S}_1 = \{s _1 \}, \mathbb{S}_2 = \{s_2 \},...,
  \mathbb{S}_n = \{s_n \} $ we define 
 \begin{equation*}
   \mathbb{S}_1 \curlyvee \mathbb{S}_2 \curlyvee ... \curlyvee \mathbb{S}_n =
   (s_1, s_2,..., s_n).
 \end{equation*}
Note that $\zeta _{2n} = (\zeta _1 \setminus \beta) \curlyvee (\zeta _3 \setminus \zeta _2) \curlyvee
... (\zeta _{2n-1} \setminus \zeta _{2n-2})$ if 
$(\zeta_1,...,\zeta _{2n}) \in \Xi$.

From the definition of the Lebesgue Poisson measure we have
 \begin{equation}\label{straitjacket}
  l_\d ^n(sym ^{-1} A) = n! \lambda (A),
 \end{equation}
 where $l_\d ^n$ is 
 the Lebesgue measure on $(\R ^\d)^{n}$.

 Define a measure $\sigma$ on 
 $\big(\prod\limits _{k=1} ^n B_{\frac r4}(x_k), \mathscr{B}( \prod\limits _{k=1} ^n B_{\frac r4}(x_k)) \big)$ by
 
 \begin{equation*}
  \begin{gathered}
     \sigma(D) =  \int\limits _{\zeta_1, \zeta _2, ... , \zeta _{2n}} 
 {Q(\beta, d \zeta _1)}   
 {Q(\zeta _1, d \zeta _2)}
 {Q(\zeta _2, d \zeta _3)} \times ...\times
{Q(\zeta _{2n-1}, d \zeta _{2n})} 
\\
 \times I_{\{(\zeta_1,...,\zeta _{2n}) \in \Xi  \}} 
I_{\{ (\zeta _2 \setminus \zeta _1) \curlyvee (\zeta _4 \setminus \zeta _3) \curlyvee
... \curlyvee (\zeta _{2n} \setminus \zeta _{2n-1}) \in  D \}}, 
\ \ \  D \in \mathscr{B} \big(\prod\limits _{k=1} ^n B_{\frac r4}(x_k) \big) .
   \end{gathered}
 \end{equation*}

We can rewrite \eqref{infidelity} as

\begin{equation}\label{rowdy}
 Q ^{2n}(\beta, A) \geq \sigma (sym ^{-1} A).
\end{equation}
 We will show that 
 \begin{equation}\label{procure}
  \sigma (D) \geq \tilde c_3 l _\d ^n (D) ,  \ \ \ D \in \mathscr{B}(\prod\limits _{k=1} ^n B_{\frac r4}(x_k))
 \end{equation}
for some constant $\tilde c_3 >0$.

The statement of the lemma is a consequence of \eqref{straitjacket}, \eqref{rowdy}
and \eqref{procure}.
To establish \eqref{procure} we only need to consider sets of the form $D_1 \times ... \times D_n$, 
$D_j \in \mathscr{B} (B_{\frac r4}(x_j))$. Define

\begin{align*}
 \Xi_{(D_1,...,D_n)}  = \bigg\{ (\zeta_1,...,\zeta_{2n}) \mid 
  \zeta _{2k-1} = \{y_1,...,y_k,x_k,...,x_n \},
  \zeta _{2k} = & \{y_1,...,y_k,x_{k+1},...,x_n \}, k=1,...,n, \\
  \text{ for some distinct } & y_k \in D_k
  \bigg\}.
\end{align*}
We have

\[
\sigma (D_1 \times ... \times D_n) = \int\limits _{\zeta_1, \zeta _2, ... , \zeta _{2n}} 
 {Q(\beta, d \zeta _1)}   
 {Q(\zeta _1, d \zeta _2)}
 {Q(\zeta _2, d \zeta _3)} \times ...
 \]
 \[
 \times
{Q(\zeta _{2n-1}, d \zeta _{2n})} 
I_{\{(\zeta_1,...,\zeta _{2n}) \in \Xi_{(D_1,...,D_n)}  \}}.
\]

Fix $z_j \in D _j$. 
Using our assumptions on $b$ and $d$, we see that
\[
Q\bigg(\{z_1,...,z_k,x_k,...,x_n \} , \big\{ \{z_1,...,z_k,x_{k+1},...,x_n \} \big\} \bigg) = 
\frac{d(x_k, \{z_1,...,z_k,x_k,...,x_n \})}{(B+D)\{z_1,...,z_k,x_k,...,x_n \}}
\]
\[
 \geq \frac{d(x_k, \{z_1,...,z_k,x_k,...,x_n \})}{\sup \{(B+D)(\eta)\mid |\eta| \leq n+1\}} \geq 
 \frac{ \inf\limits _{\eta \in \Gamma _0 (\R^\d), x \in \eta} d(x,\eta)}{\sup \{(B+D)(\eta)\mid |\eta| \leq n+1\}}\;,
\]
and 
\[
Q \bigg( \{ z_1,...,z_k,x_{k+1},...,x_n \} , \big\{ \{z_1,...,z_k,y_{k+1},x_{k+1},...,x_n \} 
\mid y_{k+1} \in D_{k+1}  \big\} \bigg) = 
\]
\[
 = \frac{\int\limits _{y \in D_{k+1}} b(y, \{z_1,...,z_k,x_{k+1},...,x_n \}) dy}
 {(B+D)(\{z_1,...,z_k,x_{k+1},...,x_n \})} \geq 
 \frac{c_3 l_\d (D_{k+1})}
 {(B+D)(\{z_1,...,z_k,x_{k+1},...,x_n \})},
\]
where $l_\d$ is the Lebesgue measure on $\R ^\d$.
Hence
\[
\sigma (D_1 \times ... \times D_n) \geq 
\bigl(
\frac{ \inf\limits _{\eta \in \Gamma _0 (\R^\d), x \in \eta} d(x,\eta)}{\sup \{(B+D)(\eta) :  |\eta| \leq n+1 \}}
\bigr) ^n
 \prod\limits _{j=1} ^n \frac{c_3 l_\d (D_{j})}
 {\sup \{(B+D)(\eta) :  |\eta| \leq n+1 \}} .
 \]

 It remains to note that 
 $\prod\limits _{j=1} ^n l_\d (D_j)  = l _\d ^n(D_1 \times ... \times D_n) $.

 \textbf{Proof of Theorem} \ref{irreducibility thm}. We will first establish $\phi$-irreducibility.
 Starting from any configuration, the process may go extinct in finite time: for all
  $\eta \in \Gamma _0 (\R^\d)$
\[
 Q _{\eta} \{ \xi _k = \varnothing \text{ for some } k>0 \} >0 .
\]
 Therefore, it is sufficient to show that
 \begin{equation}\label{tentacle}
   L(\varnothing , A) >0 \ \ \text{ whenever } \lambda (A) >0,
 A \in \mathscr{B}(\Go). 
 \end{equation}
 
Let us take $A \in \mathscr{B}(\Go)$ with $\lambda(A)>0$.
There exists $n \in \N$ and 
$\beta' \in \Gamma _0 ^{(n)}$ such that

 \begin{equation}
 \lambda (A \cap \B _{\rho}(\beta', \frac {r}{4}) ) > 0.
\end{equation}
 
 By Lemma \ref{5lem11} there exists a path from $\varnothing$ 
 to $\beta'$. Denote by $m$ the length of this path. 
 Applying Lemma \ref{5lem12} and Lemma \ref{5lem13} we get
 
 \[
  Q^{m+2n}(\varnothing, A ) \geq \int\limits _{\beta \in \B _{\rho}(\beta', \frac {r}{4})}
  Q^{m}(\varnothing, d \beta ) Q^{2n}(\beta, A) >0,
 \]
which proves \eqref{tentacle}.

 Now let us prove that $\lambda$ is a maximal irreducibility measure for 
 $(\xi _n)_{n\in \N}$. Taking into account Proposition \ref{lush}, we see that
 it suffices to show that for all $A \subset \Gamma _0 (\R^\d)$ with
 $\lambda (A) = 0$ we have

 \begin{equation} \label{ghj}
 \lambda  \{ \eta : Q(\eta , A)  >0 \} = 0.
 \end{equation}

 With no loss of generality, we assume that
 $A \subset \Gamma _0 ^{(n)} (\R^\d)$, $n\geq 2$. 
 We have $sym ^{-1} (A) \subset (\R^\d)^n$
 and $l_{\d} ^n (sym ^{-1} (A) ) = 0$. Now,
 $\eta \in \Gamma _0 ^{(n+1)} (\R^\d)$ and $Q(\eta , A) >0$ 
 if and only if $\eta$ may be represented as
 $\xi \cup \{ x \}$, 
 where $\xi \in A$, $x \in \R^\d \setminus \xi $. 
 Then we also have for any
 $y = (y_1,...,y_{n+1}) \in sym ^{-1} (\eta)$
 \[
 \check{\Pi} _j y \in sym ^{-1} (A) 
 \]
 for 
 some $j \in \{1,2,...,n+1  \}$, where 
 $\check{\Pi} _j y = (y_1,...,y_{i-1}, y_{i+1},..., y_{n+1}
 ) \in (\R^\d)^n$. 

Since
 $l _{\d} ^n (sym ^{-1} (A) ) = 0$, we also have
 $l _{\d} ^{n+1} (\check{\Pi}(\cdot) _j ^{-1} (sym ^{-1} (A)) ) = 0$,
 and consequently 
 
 \begin{equation}\label{adamant}
 \lambda  \{ \eta : \eta \in \Gamma _0 ^{(n+1)},  Q(\eta , A)  >0 \} = 0.
\end{equation}

 Similarly, if $\eta \in \Gamma _0 ^{(n-1)} (\R^\d)$ and $Q(\eta , A) >0$, then for $y \in sym ^{-1} (\eta)$

\begin{equation}\label{ascribe}
 l _{d} \{ z \in \R ^\d:  (z,y) \in  sym ^{-1} (A)  \} >0.
\end{equation}
 because a ``newly born''  point has an absolutely continuous distribution with respect to the
 Lebesgue measure on $\R ^\d$, in the sense that
 $Q(\eta, \{\eta \cup z \mid z \in D \}) = 0$ if $l_\d(D)=0$.
 However, the set of all $y$ satisfying \eqref{ascribe} has zero 
 Lebesgue measure, otherwise we would have 
 \[
 l _{\d} ^n ( sym ^{-1} (A)) = 
 \int l_{\d} ^{n-1}(dy)l_{\d}\{z: (z,y) \in sym^{-1}(A) \} >0.
 \]
 Therefore, 
\begin{equation}\label{serene}
 \lambda  \{ \eta : \eta \in \Gamma _0 ^{(n-1)},  Q(\eta , A)  >0 \} = 0.
\end{equation}
Note
that in cases $n=0,1$ some changes should
be made in the proofs of \eqref{adamant}, \eqref{serene}, 
because of 
the special structure of $\Gamma _0 ^{(0)} (\R^\d) = 
\{ \varnothing \}$. 

Since 
$$  \{ \eta : \eta \in \Gamma _0 ^{(k)},  P(\eta , A)  >0 \} = \varnothing,$$
$k \ne n-1, n+1$, $n \geq 0$,
 \eqref{adamant} and \eqref{serene} imply \eqref{ghj}.
$\qed$

  \section*{Acknowledgement}

  Viktor Bezborodov is  supported by
the Department of Computer Science at the University of Verona.
The authors are grateful to the anonymous referee(s) for their constructive comments.
 
\bibliographystyle{alpha}
\bibliography{SinusUpdated}

\end{document}